\documentclass[11pt]{article}

\usepackage[a4paper,margin=24mm]{geometry}
\usepackage{amsmath,amssymb,amsthm,mathtools}
\usepackage{subfiles}
\usepackage{xcolor}
\usepackage{microtype}
\usepackage[pdfusetitle,colorlinks=true,linkcolor=blue,citecolor=blue,
urlcolor=blue]{hyperref}
\usepackage{interval}
\intervalconfig{soft open fences}
\usepackage{derivative}
\usepackage{mleftright}
\usepackage{tikz}
\usepackage{graphicx}
\usetikzlibrary{arrows.meta,positioning,calc}

\newtheorem{theorem}{Theorem}[section]
\newtheorem{lemma}[theorem]{Lemma}
\newtheorem{proposition}[theorem]{Proposition}

\theoremstyle{remark}
\newtheorem{remark}[theorem]{Remark}

\newcommand{\E}{\mathbb E}
\newcommand{\Var}{\operatorname{Var}}
\newcommand{\tr}{\operatorname{tr}}
\newcommand{\I}{\mathbb I}
\newcommand{\N}{\mathcal N}
\newcommand{\calI}{\mathcal I}
\newcommand{\norm}[1]{\left\lVert #1\right\rVert}

\title{A Central Limit Theorem for High-Temperature SYK Model: A case study on mathematical research using AI and Lean}
\author{Yoshinori Kamijima\thanks{Information Networking for Innovation and Design, Toyo University. E-mail: kamijima@iniad.org}
\and Shuta Nakajima\thanks{Department of Mathematics, Keio University. E-mail: njima@keio.jp}
\and Ryoichiro Noda\thanks{Faculty of Science and Engineering, Waseda University.
E-mail: r\textunderscore noda@aoni.waseda.jp}
\and K\^{o}hei Sasaya\thanks{Graduate School of Mathematical Sciences, The University of Tokyo. E-mail: sasaya@g.ecc.u-tokyo.ac.jp}
}
\date{}

\begin{document}

\maketitle

\begin{abstract}
The Sachdev--Ye--Kitaev (SYK) model is a quantum many-body model of Majorana
fermions with all-to-all random interactions. It has attracted wide interest in condensed matter and high-energy physics as a solvable model of non-Fermi liquids and holography. Fix an even integer $q\geq4$, and let $Z_N(\beta)$ be the
partition function of the Gaussian $q$-body SYK model.  We prove that for all sufficiently small inverse temperatures
$\beta > 0$, there exists $\rho=\rho(\beta)>0$ such that
\[
N^{q/2 -1}(\log Z_N(\beta)-\E[\log Z_N(\beta)])
  \Longrightarrow \N(0,\rho^2).
\]
In addition to the mathematical results, we discuss how the proofs were developed and refined through interaction between human researchers and AI, and we describe the subsequent formalization of the main arguments in Lean 4.
This work thus serves as a case study of a research workflow that integrates human insight, AI-assisted exploration, iterative proof refinement, and formal verification.
\end{abstract}



\section*{Research process with AI and Lean}

In the past year, AI has been used more and more in mathematical research, often together with the Lean proof assistant \cite{GGTW2025,OpenAI}. 
Lean is a programming language and proof assistant that verifies the correctness of formalized proofs.
We had been studying the Sachdev--Ye--Kitaev (SYK) model together in a seminar, during which we became  interested in establishing  a central limit theorem (CLT) for its free energy. 
We decided to investigate this problem as a case study of a mathematical research workflow involving AI and Lean.

We first attempted to prove the CLT for the free energy directly through discussions with AI (ChatGPT 5.5 Pro), but these initial attempts did not lead to a proof.
We then proposed studying the CLT for the partition function itself, rather than directly for the free energy.
After suggesting this alternative approach to AI and continuing the discussion, we obtained proofs of the CLTs for the partition function and, consequently, the free energy. 
We subsequently examined the AI-generated proof carefully and worked on verifying and refining it.
In particular, the initial proof of a key estimate relied on a complicated analytic argument. 
Through further discussions with AI (ChatGPT 5.6 Pro), we finally developed a simpler combinatorial argument, which is the approach adopted in Section~\ref{sec:Rademacher-remainder} to control the higher-order chaos terms.
We then reorganized the paper so that the proofs are more accessible to human readers. In parallel, we formalized the refined proofs in Lean: the AI generated the code, and we reviewed it.

We now describe how we used AI to develop the Lean formalization. 
We first asked AI to formalize the model and the statements of our main results, and checked that these accurately represented the mathematics we intended. 
We then used AI to write the corresponding proofs in Lean. 
To guide this work, we prepared written instructions and supporting documents specifying how AI should proceed. 
In particular, these documents were intended to reduce the risk of AI changing the statements to be proved or introducing additional assumptions. 
Although preparing these documents requires some initial effort, it can reduce the human effort needed to identify and correct such problems later.
Finally, we used AI to reorganize and simplify the Lean code to make it easier to read, while preserving its mathematical content.
For further details, see \url{https://github.com/TokyoProbabilityGroup/SYK_CLT}.

\begin{figure}[t]
    \centering
    \resizebox{0.9\linewidth}{!}{%
    \begin{tikzpicture}[
        node distance=10mm and 18mm,
        >=Latex,
        human/.style={
            draw=red,
            thick,
            rounded corners=2pt,
            rectangle,
            align=center,
            minimum height=10mm,
            inner sep=6pt,
            fill=red!10
        },
        ai/.style={
            draw=blue,
            thick,
            rounded corners=2pt,
            rectangle,
            align=center,
            minimum height=10mm,
            inner sep=6pt,
            fill=blue!10
        },
        final/.style={
            draw=black!70,
            thick,
            rounded corners=2pt,
            rectangle,
            align=center,
            minimum height=10mm,
            inner sep=6pt,
            fill=white
        },
        arr/.style={
            -{Latex[length=2.2mm]},
            shorten <=1.5mm,
            shorten >=1.5mm,
            thick
        },
        darr/.style={
            <->,
            shorten <=1.5mm,
            shorten >=1.5mm,
            thick
        }
    ]

    \node[human] (study) {Study the research topic};
    \node[ai, below=of study] (problem) {Formulate the problem and explore solutions};
    \node[human, below=of problem] (polish) {Refine and polish the proof};
    \node[ai, below=of polish] (lean) {Generate Lean code from the polished proof};

    \node[human] (check) at ($(lean.south)+(-45mm,-23mm)$) {Examine the Lean code};
    \node[human] (manuscript) at ($(lean.south)+(45mm,-23mm)$) {Refine the manuscript for human readers};

    \coordinate (mid) at ($(check)!0.5!(manuscript)$);
    \node[final, below=18mm of mid] (complete) {Completion of the work};

    \draw[arr] (study) -- (problem);
    \draw[arr] (problem) -- (polish);
    \draw[arr] (polish) -- (lean);

    \draw[arr] (problem.north east)
        .. controls +(40mm,20mm) and +(40mm,-20mm) .. (problem.south east)
        node[pos=0.5,human] (feedback) {Human feedback};

    \draw[arr, shorten <=4.5mm, shorten >=4.5mm] (lean.south) -- (check.north);
    \draw[arr, shorten <=4.5mm, shorten >=4.5mm] (lean.south) -- (manuscript.north);

    \draw[darr, dashed] (check.east) -- node[above] {Mutual feedback} (manuscript.west);

    \draw[arr, shorten <=4.5mm, shorten >=4.5mm] (check.south) -- (complete.north);
    \draw[arr, shorten <=4.5mm, shorten >=4.5mm] (manuscript.south) -- (complete.north);

    \end{tikzpicture}
    }
    \caption{
        A schematic workflow for mathematical research with AI and Lean.
Boxes shaded in red indicate steps carried out primarily by human researchers,
whereas boxes shaded in blue indicate steps in which AI may play an essential role.
    }
    \label{fig:workflow}
\end{figure}
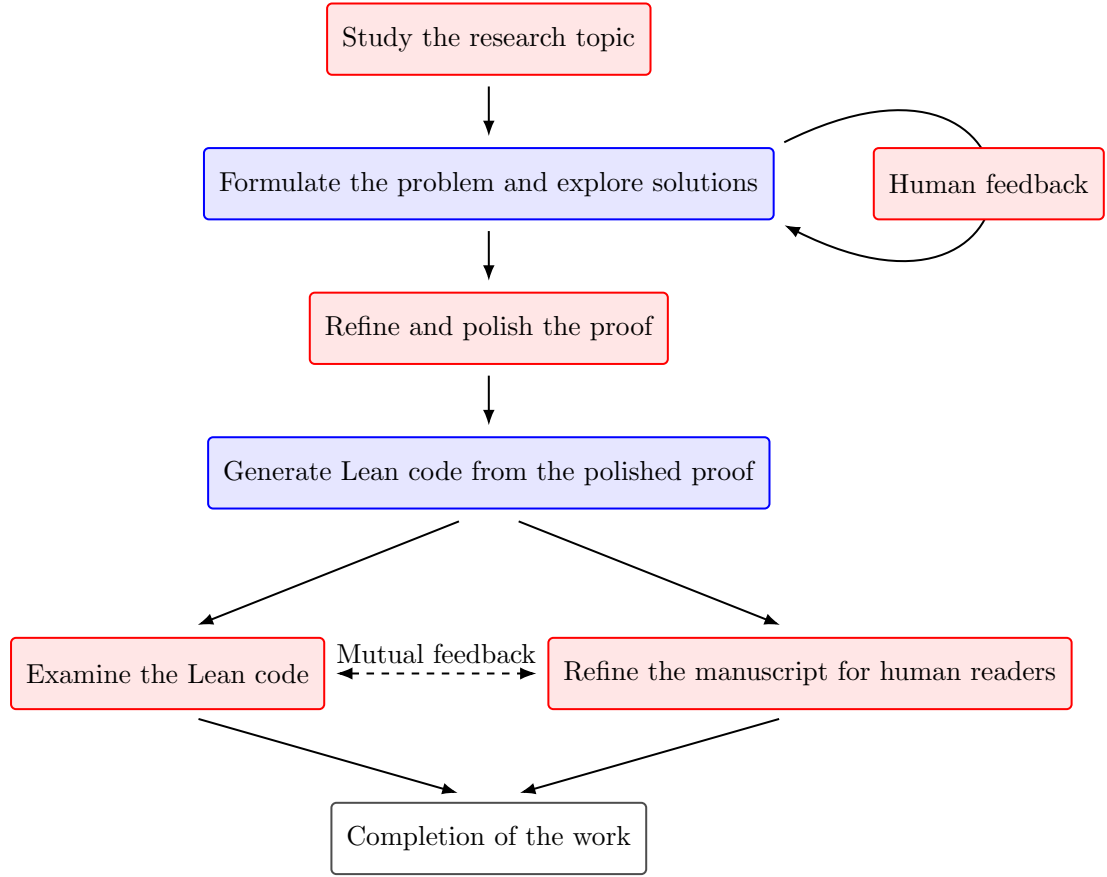

We expect that the workflow illustrated in Figure~\ref{fig:workflow} will be useful particularly for master's and PhD students beginning mathematical research with AI. Moreover, the process of examining and discussing AI-generated arguments may itself help researchers deepen their mathematical understanding.  
From a practical perspective, refining a mathematical proof before formalizing it in Lean can help produce simpler and better-structured code.

\section{Introduction and main results}
\label{sec:introduction}

\subsection{The SYK model}\label{subsec:SYK-model}
The Sachdev--Ye--Kitaev (SYK) model is a quantum-mechanical model of \(N\) Majorana fermions with random many-body interactions. It originates in the work of Sachdev and Ye on random quantum magnets \cite{SachdevYe1993} and in Kitaev’s subsequent formulation in terms of Majorana fermions \cite{Kitaev2015}. Later work by Maldacena, Stanford, and others clarified its large-\(N\) behavior and revealed deep connections to quantum chaos and two-dimensional quantum gravity \cite{MaldacenaStanford2016}. These features make the SYK model a useful setting for studying strongly interacting quantum many-body systems.

Depending on how the inverse temperature \(\beta=\beta_N\) scales with \(N\), the model exhibits several distinct regimes. For fixed \(\beta\), one obtains the standard finite-temperature large-\(N\) regime; when \(1\ll \beta\ll N\), the system enters the conformal regime, characterized by strong interactions and non-Fermi liquid behavior. At lower temperatures, \(N\ll \beta\), the quantum-Schwarzian regime emerges, in which finite-\(N\) effects and quantum fluctuations of the Schwarzian mode play a central role. Finally, when \(\beta\) is exponentially large in \(N\), the partition function is mainly determined by the eigenvalues near the edge of the spectrum. In this regime, random matrix theory  describes their fluctuations. Hence, the SYK model shows several temperature crossovers and allows analytic calculations in the scaling limits.

Previous mathematical work has studied the limiting free energy of the SYK model for fixed $\beta$ \cite{GPSZ2026,He2026}. Moreover, Anschuetz, Gamarnik, and Kiani \cite{AGK2025} established a universality result showing that the limiting free energy is the same for a broad class of interaction distributions.  This naturally raises the question of how physical quantities depend on the interactions. Motivated by this question, we study fluctuations of the free energy and show that, for sufficiently small \(|\beta|\), the suitably centered and rescaled free energy converges to a Gaussian distribution.

Fix even integers $q,N$ such that $4\leq q\leq N$. Set $[N] \coloneqq \{1,2,\ldots,N\}$ and write
\[
  \calI_N \coloneqq \{I\subset[N]:|I|=q\},\qquad
  \sigma_N^2 \coloneqq \frac{(q-1)!}{N^{q-1}},
\]
where $|I|$ denotes the cardinality of the set $I$.
Let $G \coloneqq (g_I)_{I\in\calI_N}$ be a family of independent standard Gaussian random
variables. We denote the corresponding probability measure and expectation by $\mathbb P$ and $\E$, respectively. Following \cite[Section~2]{GPSZ2026}, we let $\psi_1,\ldots,\psi_N$ be Majorana matrices of size $2^{N/2} \times 2^{N/2}$, which satisfy the following anticommutation relation:
\begin{equation}
  \label{eq: clifford relation}
  \{\psi_i,\psi_j\} \coloneqq \psi_i \psi_j + \psi_j \psi_i =2\delta_{ij}\I,
\end{equation}
where $\I$ denotes the identity matrix.
Note that the Majorana matrices considered here are $\sqrt{2}$ times those used in \cite{GPSZ2026,He2026}; this is merely a difference in normalization. 
For
$I = \{i_1, \ldots, i_q\} \in \calI_N$ with $i_1 < \cdots < i_q$, set
\[
  \psi_I \coloneqq \mathrm i^{q/2}\psi_{i_1}\cdots\psi_{i_q},
\]
where $\mathrm i$ denotes the imaginary unit. Since Majorana matrices are Hermitian, it follows from \eqref{eq: clifford relation} that 
$\psi_I$ is also Hermitian. Moreover, $\psi_I^2=\I$ and, in particular, $\norm{\psi_I}_\infty=1$, where $\norm{\,\cdot\,}_\infty$ denotes the operator norm of a matrix.  
The Hamiltonian of the SYK model is defined by 
\[
  H_N \coloneqq H_N(G) \coloneqq \sigma_N\sum_{I\in\calI_N}g_I\psi_I.
\]
The corresponding partition function is then given by
\[
  Z_N(\beta) \coloneqq Z_N(\beta;G) \coloneqq \tr e^{-\beta H_N(G)},\qquad \beta \in \mathbb R.
\]
Here, $\tr \coloneqq 2^{-N/2}\operatorname{Tr}$ denotes the normalized trace; in
particular, $Z_N(0;G)=1$.  

\subsection{Main results}
Here, we state our main results. 
We first recall the convergence of the annealed free energy, which was established for sufficiently small $|\beta|$ in \cite{GPSZ2026} and recently extended to all $\beta$ in \cite{He2026}. Moreover, in \cite{He2026}, the limit is identified with the Schwinger--Dyson pressure. Throughout this article, all limits as $N \to \infty$ are taken through even integers.
\begin{theorem}[{\cite[Theorem~1.1 and Proposition~3.3]{He2026}}]\label{thm: input}
For any $\beta\in \mathbb R$, there exists  $\ell_q(\beta)$ such that
\begin{align*}
    \lim_{N\to \infty} N^{-1} \log \E[Z_N] =\ell_q(\beta).
\end{align*}
Moreover, $\ell_q(\beta)$ is differentiable with respect to $\beta$ on $\mathbb R$.
\end{theorem}

Our first main result concerns the fluctuations of the partition function ${Z_N(\beta)}$. For $a > 0$, we write $\N(0, a)$ for the centered Gaussian distribution with variance $a$. We denote by $\Longrightarrow$ convergence in distribution.
\begin{theorem}
\label{thm:partition-clt}
There exists $\beta_q\in (0,\infty)$ such that for any $\beta\in (-\beta_q,\beta_q) \setminus \{0\}$,  
 \begin{equation}\label{eq:scaled-relative-clt}
  N^{q/2-1}\left( \frac{Z_N(\beta)}{\E[Z_N(\beta)]} - 1\right)
  \Longrightarrow
  \,\N\left(0,\frac{q!}{2}\beta^2(\ell_q'(\beta))^2\right),
\end{equation}
where  $\beta \ell_q'(\beta)>0$ for any $\beta\in \mathbb R\setminus \{0\}$; in particular, the limiting distribution is non-degenerate. Moreover,
\begin{equation}\label{eq:standardized-clt}
  \frac{Z_N(\beta)-\E[Z_N(\beta)]}
  {\sqrt{\Var Z_N(\beta)}}
  \Longrightarrow\N(0,1),
\end{equation}
where $\Var Z_N(\beta)$ denotes the variance of $Z_N(\beta)$.
\end{theorem}

Combining the above theorem with the asymptotic relation $\log x \sim x-1$ as $x \to 1$, we obtain the following result on the fluctuations of the free energy.
\begin{theorem}\label{thm:free-energy-clt}
With the same $\beta_q$ as in Theorem~\ref{thm:partition-clt}, for any $\beta\in (-\beta_q,\beta_q) \setminus \{0\}$,
\[
  N^{q/2-1}\left(
    \log Z_N(\beta)-\E[\log Z_N(\beta)]
  \right)
  \Longrightarrow
  \,\N\left(0,\frac{q!}{2}\,\beta^2(\ell_q'(\beta))^2\right).
\]
\end{theorem}

\begin{remark}
The proof of Theorem~\ref{thm:partition-clt} is based on a polynomial chaos expansion of $Z_N(\beta)$, and the smallness of $|\beta|$ is used to control the higher-order chaos terms; see Proposition~\ref{prop:higher-chaos-closure} below.
We do not pursue an extension to all $\beta\in\mathbb R$ here, nor do we make a conjecture as to whether the same fluctuation results hold at all temperatures.
\end{remark}

\begin{remark}
In related work, Feng, Tian, and Wei \cite{FTW2018} proved a central limit theorem for linear eigenvalue statistics of the Hamiltonian $H_N$. Their result applies to a fixed test function \(f\), under suitable regularity assumptions on \(f\) and \(f'\). In their normalization, however, the partition function corresponds to
\[
f_{N,\beta}(x) \coloneqq e^{-\beta \sqrt{N} x},
\]
so both the test function and its derivative depend on \(N\). Therefore, their result does not apply directly to \(Z_N(\beta)\).
\end{remark}

\subsection{Outline of the proofs and organization of the paper}

We give an outline of the proofs of the main results.
Define the normalized partition function by
\begin{equation} \label{eq:Y_N}
  Y_N(\beta) \coloneqq Y_N(\beta; G) \coloneqq {Z_N(\beta)}/{\E[Z_N(\beta)]}.
\end{equation}
By applying a polynomial chaos expansion to $Y_N(\beta)$, we write 
\[
  Y_N(\beta) = \sum_{k = 0}^\infty P_kY_N(\beta),
\]
where $P_k Y_N(\beta)$ denotes the $k$-th chaos component.
(The precise definition of the above expansion is given in Section~\ref{sec:chaos-decomposition}.) The normalization gives $P_0 Y_N(\beta) = \E [Y_N(\beta)] = 1$. Moreover, in {Proposition~\ref{prop:exact-second-chaos}}, using the anticommutation relation~\eqref{eq: clifford relation} and Gaussian integration by parts, we deduce that $P_1 Y_N(\beta) = 0$ and
\[ 
  P_2Y_N(\beta) = A_N(\beta) Q_N,
\]
where
\begin{equation}\label{eq:finite-N-FTW-coefficient}
      A_N(\beta) \coloneqq \frac{\beta}{2}\partial_\beta\log \E[Z_N(\beta)],
  \qquad
    Q_N \coloneqq \frac1{|\mathcal I_N|}\sum_{I\in\mathcal I_N}(g_I^2-1).
\end{equation}
Hence, we obtain
\begin{equation}\label{eq:intro-chaos-expansion}
  Y_N(\beta)-1=
  A_N(\beta)Q_N+\mathcal R_N(\beta),
\end{equation}
where $\mathcal R_N(\beta) \coloneqq \sum_{k \geq 3} P_k Y_N(\beta)$.
Theorem~\ref{thm: input} and the convexity of $\log \E[Z_N(\beta)]$ in $\beta$ imply that
\begin{equation}\label{eq:intro-A_N-limit}
  \frac{A_N(\beta)}N
  = 
  \frac{\beta}{2}\cdot \frac{1}{N}\partial_\beta\log \E[Z_N(\beta)]
  \xlongrightarrow{N \to \infty} \frac{\beta}{2}\ell_q'(\beta),
\end{equation}
as shown in Proposition~\ref{prop:annealed-free-energy}.
Since the usual central limit theorem applies to $Q_N$, Theorem~\ref{thm:partition-clt} follows from \eqref{eq:intro-chaos-expansion} and \eqref{eq:intro-A_N-limit}, combined with  
\begin{equation}\label{eq:intro-higher-chaoses}
  \E \bigl[ \mathcal R_N(\beta)^2 \bigr]= o(N^{-(q-2)}),\qquad \text{as } N \to \infty,
\end{equation}
where $o(\cdot)$ denotes the usual Landau notation.
The above estimate is proved in Proposition~\ref{prop:higher-chaos-closure} 
 {by using a combinatorial argument to bound the number of index configurations that give a non-zero contribution to the higher-order chaoses.}

To obtain Theorem~\ref{thm:free-energy-clt}, observe that
\begin{align*}
  &N^{q/2-1}\bigl(\log Z_N(\beta)-\E[\log Z_N(\beta)]\bigr)\\
  &=N^{q/2-1}\bigl(\log Y_N(\beta)-(Y_N(\beta)-1)\bigr)
   +N^{q/2-1}(Y_N(\beta)-1)-N^{q/2-1}\E[\log Y_N(\beta)].
\end{align*}
The Taylor expansion $\log(1+x)-x=O(x^2)$ as $x \to 0$, together with Theorem~\ref{thm:partition-clt}, shows that the first term  on the right-hand side converges to zero in probability. Moreover, using the estimates obtained in the proof of \eqref{eq:intro-higher-chaoses}, we deduce that $\E[(Y_N(\beta)-1)^2]=O(N^{2-q})$ and
$\mathbb P(Y_N(\beta)<1/2)=O(N^{-3(q/2-1)})$. Together with
$Z_N(\beta)\geq1$, these estimates give
\[
  \E[\log Y_N(\beta)]=O\!\left(N^{2-q}\right)
  =o\!\left(N^{1-q/2}\right).
\]
This completes the proof of Theorem~\ref{thm:free-energy-clt}; see Section~\ref{sec:main-results} for details.

{The rest of the paper is organized as follows. In
Section~\ref{sec:chaos-decomposition}, we introduce the polynomial chaos expansion
of $Z_N(\beta)/\mathbb E Z_N(\beta)$ and compute its first two chaos components. We show that the first
chaos vanishes and that the second chaos equals $A_N(\beta)Q_N$, where $A_N(\beta)$ and $Q_N$ 
are defined in~\eqref{eq:finite-N-FTW-coefficient}. In
Section~\ref{sec:Rademacher-remainder}, we show that the chaos components of order
three and higher are negligible compared with the second chaos. In Proposition~\ref{prop:annealed-free-energy} of Section~\ref{sec:main-results},
we derive the asymptotics of $A_N(\beta)$ from the annealed free energy
estimates of~\cite{He2026} and show that $\beta\ell_q'(\beta)>0$ for all
$\beta\neq 0$.  Combining these results with the classical central limit theorem for $Q_N$, we prove
Theorems~\ref{thm:partition-clt} and~\ref{thm:free-energy-clt}.}

\section{Polynomial chaos expansion}
\label{sec:chaos-decomposition}
We call $\alpha=(\alpha_I)_{I\in\mathcal I_N} \in (\mathbb N\cup\{0\})^{\mathcal I_N}$  a multi-index. We define the order of $\alpha$ as $|\alpha| := \sum_{I\in\mathcal I_N} \alpha_I $.
  {For $k\in\mathbb Z_{\geq0}$, the probabilists' Hermite
 polynomial of degree $k$ is defined by
 \[
   {\rm He}_k(x):=(-1)^k e^{x^2/2}\frac{d^k}{dx^k}e^{-x^2/2}.
 \]
 These polynomials are orthogonal with respect to the standard Gaussian
 measure; in particular, ${\rm He}_0(x)=1$, ${\rm He}_1(x)=x$, and
 ${\rm He}_2(x)=x^2-1$.  Their multivariate products satisfy}
\begin{equation*}
  H_\alpha(G):=\prod_{I\in\calI_N}{\rm He}_{\alpha_I}(g_I),
  \qquad
  \E[H_\alpha H_\eta]
  =\alpha!\,\mathbf 1_{\{\alpha=\eta\}},
\end{equation*}
where $\alpha!:=\prod_I\alpha_I!${ and
$\mathbf 1_{\{\alpha=\eta\}}$ is the indicator function}. Since the $H_\alpha$ form a complete
orthogonal basis of $L^2(\gamma_{\mathcal I_N})$, where
$\gamma_{\mathcal I_N}$ denotes {the} standard Gaussian measure on
$\mathbb R^{\mathcal I_N}$, we have
\begin{align}\label{eq: chaos}
Z_N(\beta;G) ={\sum_{k=0}^\infty \sum_{\alpha:|\alpha|=k} \frac{\E[Z_N(\beta;G)H_\alpha(G)]}{\alpha!} H_\alpha(G)},
\end{align}
where the second sum runs over all multi-indices $\alpha$ such that $\sum_{I} \alpha_I = k$.  Let $P_r$ denote the orthogonal projection in
$L^2(\gamma_{\mathcal I_N})$ onto the $r$-th  chaos {component}.
Then, we can write
\[
  P_r Z_N(\beta;G)
  =\sum_{|\alpha|=r}
  	\frac{\E[Z_N(\beta;G)H_\alpha(G)]}{\alpha!}H_\alpha(G).
\]
The aim of this section is to prove the following proposition for the normalized partition function.
\begin{proposition}\label{prop:exact-second-chaos}
	For every $\beta\in \mathbb R$,
	\[
	P_1Y_N(\beta)=0\quad \text{and} \quad P_2 Y_N(\beta)
	=A_N(\beta)Q_N,
	\]
	where $A_N(\beta)$ and $Q_N$ are defined by \eqref{eq:finite-N-FTW-coefficient}.
\end{proposition}

{We call}  a multi-index $\alpha$ \emph{even} if $\sum_{I\colon x \in I}
\alpha_I$ is even for every $x\in[N]$, and odd otherwise. We write $\mathfrak S_N$ for the symmetric group of degree $N$. 
\begin{lemma}\label{lem: parity rule lemma}
    For $s=(s_1,\ldots,s_N)\in\{\pm1\}^N$, define
$\chi_s(I):=\prod_{x\in I}s_x$.
 {For $\pi\in\mathfrak S_N$ and $I=\{i_1,\ldots,i_q\}\in\calI_N$ with $i_1<\cdots<i_q$, we also
define $\epsilon_\pi(I)\in\{\pm1\}$ by
\[
  \mathrm i^{q/2}\psi_{\pi(i_1)}\cdots\psi_{\pi(i_q)}
  =\epsilon_\pi(I)\psi_{\pi(I)},
\]
whose existence follows from the anticommutation relation \eqref{eq: clifford relation}.
Set
\[
  T_{\pi,s}G
  :=\bigl(\chi_s(I)\epsilon_\pi(I)g_{\pi(I)}\bigr)_{I\in\calI_N}.
\]
Then, for every $\beta\in\mathbb R$ and $G\in\mathbb R^{\calI_N}$,
\begin{equation}\label{eq:signed-permutation-invariance}
  Z_N(\beta;T_{\pi,s}G)=Z_N(\beta;G).
\end{equation}
}
    In particular, if $\alpha$ is odd, then
   \begin{equation}\label{eq:odd-coefficient-zero}
  \E[Z_N(\beta;G)H_\alpha(G)]=0.
   \end{equation}
    
\end{lemma}
\begin{proof}
We define 
\[
  \Gamma
  := 
  \psi_1\cdots\psi_N.
\]
Note that $\Gamma^*\Gamma$ is {the} identity {matrix}, so $\Gamma$ is unitary. Since $N$ is even, moving $\psi_j$ through the other $N-1$ Majorana
{matrices}  produces the sign $(-1)^{N-1}=-1$.  Hence, for any $j\in[N],$
\[
  \Gamma\psi_j\Gamma^*=  -\psi_j \Gamma \Gamma^* =-\psi_j.
\]
{For a real unit vector $v=(v_1,\ldots,v_N)$, set
\[
  C_v:=\sum_{j=1}^N v_j\psi_j,
  \qquad U_v:=\Gamma C_v.
\]
The anticommutation relation~\eqref{eq: clifford relation} gives
\[
  C_v^*=C_v,\qquad C_v^2=\mathbb I,\qquad
  C_v\psi_j+\psi_jC_v=2v_j\mathbb I\quad(j\in[N]).
\]
Thus $C_v$ and $U_v$ are unitary, and
\begin{equation}\label{eq:reflection-conjugation}
  U_v\psi_jU_v^*
  =\Gamma(2v_jC_v-\psi_j)\Gamma^*
  =\psi_j-2v_jC_v.
\end{equation}
In particular, taking $v=e_x$, where $e_x$ is the $x$-th
standard basis vector of $\mathbb R^N$, gives
$U_x:=\Gamma\psi_x$ and} 
\begin{align}\label{eq: conjugation rule}
      U_x\psi_jU_x^*
  =
  \begin{cases}
    -\psi_x,&j=x,\\
    \psi_j,&j\neq x.
  \end{cases}
\end{align}

Let $S_s:=\{x\in[N]:s_x=-1\}$ and let $U_s$ be the product of the unitaries
$U_x$ over $x\in S_s$:
$$U_s:= \prod_{x\in S_s} U_x,$$
where the product is in any fixed order.  Applying \eqref{eq: conjugation rule} successively yields that for any $j\in[N],$
\[
  U_s\psi_jU_s^*=s_j\psi_j.
\]
{For distinct $a,b\in[N]$, taking
$v=(e_a-e_b)/\sqrt2$ in~\eqref{eq:reflection-conjugation} gives
\begin{equation}\label{eq:unitary-transposition}
 V_{ab}:=\Gamma\frac{\psi_a-\psi_b}{\sqrt2},
\qquad
V_{ab}\psi_jV_{ab}^*
=\begin{cases}
	\psi_b,&j=a,\\
	\psi_a,&j=b,\\
	\psi_j,&j\notin\{a,b\}.
\end{cases}
\end{equation}
Writing $\pi=\tau_m\circ\cdots\circ\tau_1$ as a product of
transpositions and setting $V_\pi:=V_{\tau_m}\cdots V_{\tau_1}$,
we obtain $V_\pi\psi_jV_\pi^*=\psi_{\pi(j)}$ by \eqref{eq:unitary-transposition}.
For $\pi=\mathrm{id}$, we take $V_\pi=\mathbb I$.
The unitary $V_{\pi,s}:=V_\pi U_s$ therefore satisfies
\[
  V_{\pi,s}\psi_jV_{\pi,s}^*=s_j\psi_{\pi(j)}.
\]
Recall that $\psi_I=\mathrm i^{q/2}\psi_{i_1}\cdots\psi_{i_q}$
for $I=\{i_1,\ldots,i_q\}$ with $i_1<\cdots<i_q$.
By the definitions of $\chi_s(I)$ and $\epsilon_\pi(I)$,
\begin{align*}
  V_{\pi,s}\psi_IV_{\pi,s}^*
  &=\mathrm i^{q/2}
    (V_{\pi,s}\psi_{i_1}V_{\pi,s}^*)\cdots
    (V_{\pi,s}\psi_{i_q}V_{\pi,s}^*)\\
  &=\chi_s(I)\epsilon_\pi(I)\psi_{\pi(I)}.
\end{align*}
Since $I\mapsto\pi(I)$ is a permutation of $\calI_N$, we have
\begin{align*}
  V_{\pi,s}H_N(T_{\pi,s}G)V_{\pi,s}^*
  &=\sigma_N\sum_{I\in\calI_N}
    \chi_s(I)\epsilon_\pi(I)g_{\pi(I)}
    V_{\pi,s}\psi_IV_{\pi,s}^*\\
  &=\sigma_N\sum_{I\in\calI_N}
    \bigl(\chi_s(I)\epsilon_\pi(I)\bigr)^2
    g_{\pi(I)}\psi_{\pi(I)}
    =H_N(G).
\end{align*}
Hence, we conclude that
\[
Z_N(\beta;T_{\pi,s}G)
=\tr\!\left(
V_{\pi,s}e^{-\beta H_N(T_{\pi,s}G)}V_{\pi,s}^*
\right)
=\tr e^{-\beta H_N(G)}
=Z_N(\beta;G),
\]
because the normalized trace is invariant under unitary conjugation
and
\[
V_{\pi,s}e^{-\beta H_N(T_{\pi,s}G)}V_{\pi,s}^*
=e^{-\beta V_{\pi,s}H_N(T_{\pi,s}G)V_{\pi,s}^*}
=e^{-\beta H_N(G)}.
\]
This proves~\eqref{eq:signed-permutation-invariance}.}

Assume that $\alpha$ is odd so that  we can take $x$ such that
$\sum_{I:x\in I}\alpha_I$ is odd.  We take $s$ such that $s_y=1$ for
$y\neq x$ and $s_x=-1$.  Then,
 \[
  H_\alpha((\chi_s(I)g_I)_I)
  =
  \left(
    \prod_{x=1}^N
    s_x^{{\sum_{I:x\in I}}\alpha_I}
  \right)H_\alpha(G) = -H_\alpha(G).
\]
{
By~\eqref{eq:signed-permutation-invariance} with $\pi=\mathrm{id}$,
$Z_N(\beta;T_{\mathrm{id},s}G)=Z_N(\beta;G)$.
Since $T_{\mathrm{id},s}G$ and $G$ have the same distribution}, we have
\begin{align*}{
\E[Z_N(G) H_\alpha(G)] = \E[Z_N(\chi_s G) H_\alpha(\chi_s G)] = -
\E[Z_N(G) H_\alpha(G)]}, 
\end{align*}
which gives the second claim.
\end{proof}

For a multi-index $\alpha \in \{0,1,\ldots\}^{\mathcal I_N}$, we denote the associated multi-index partial derivative of a function $f\colon \mathbb R^{\mathcal I_N} \to \mathbb R$ as $\partial^\alpha f$.

\begin{lemma}\label{lem:Gaussian-Hermite-IBP}
Let $f\colon\mathbb R^{\mathcal I_N}\to\mathbb R$ be of class $C^{|\alpha|}$,
and suppose that $f$ and its multi-index partial derivatives up to order $\alpha$
have at most exponential growth in {$|G|$}{, where $|G|$ denotes the Euclidean norm}: for any multi-index $\eta=(\eta_I)_{I\in\mathcal I_N}$ with $\eta_I\leq \alpha_I$, there exists a  constant $C_\eta>0$ such that
\begin{equation*}
  |(\partial^\eta f)(G)| \leq {C_\eta e^{C_\eta |G|}}.
\end{equation*}
Then
\begin{equation}\label{eq:Gaussian-Hermite-IBP}
  \E\!\left[H_\alpha(G)f(G)\right]
  =\E\!\left[\partial_g^\alpha f(G)\right].
\end{equation}
In particular, \eqref{eq:Gaussian-Hermite-IBP} holds for
$f(G)=Z_N(\beta;G)$.
\end{lemma}
\begin{proof}
Recall that the one-dimensional integration-by-parts identity yields 
\[{
\E\!\left[\operatorname{He}_k(g)h(g)\right] = \E\!\left[h^{(k)}(g)\right]},
\]
for every function ${h}\in C^k(\mathbb R)$ satisfying the stated growth
condition~\cite[Equation~(2.3.15)]{S2023}. {Here, $g$ is a standard Gaussian random variable.} {By the definition of $H_\alpha(G)$, applying} this identity successively to the independent Gaussian
coordinates $(g_I)_{I\in\mathcal I_N}$ gives
\begin{align*}
\E\!\left[
H_\alpha(G)f(G)
\right]
&=
\E\!\left[
\left(
\prod_{I\in\mathcal I_N}
\operatorname{He}_{\alpha_I}(g_I)
\right)
f(G)
\right] \\
&=
\E\!\left[
\left(
\prod_{I\in\mathcal I_N}
\partial_{g_I}^{\alpha_I}
\right)
f(G)
\right] \\
&=
\E\!\left[
\partial_g^\alpha f(G)
\right].
\end{align*}
The growth assumption justifies all integrations by parts. Finally, for fixed
$N$ and $\beta$, the partition function $Z_N(\beta;G)$ is a smooth function
of finitely many Gaussian coordinates, and all of its derivatives have at
most exponential growth in  {$|G|$}.  
Thus it satisfies the assumption.
\end{proof}

\begin{proof}[{Proof { of Proposition \ref{prop:exact-second-chaos}}}]
For {any} multi-index $\alpha$ with $|\alpha|=1$, since $\alpha$ is always odd, the first claim follows from Lemma~\ref{lem: parity rule lemma}. 

{Note that, if a multi-index $\alpha$ satisfies $|\alpha|=2$ and $\alpha_I=\alpha_J=1$ for some $I,J\in\calI_N$ with $I\ne J$, then $\sum_{K: x\in K} \alpha_K=1$ for any $x\in I\setminus J$, and hence $\alpha$ is odd.
Fix any $J,J'\in\calI_N$ and choose $\pi\in\mathfrak S_N$
such that $\pi(J)=J'$.
Let $G_\pi:=T_{\pi,\mathbf 1}G$, where $\mathbf 1=(1,\ldots,1)$.
By~\eqref{eq:signed-permutation-invariance} with $s=\mathbf 1$,
$Z_N(\beta;G_\pi)=Z_N(\beta;G)$.
Since $G$ and $G_\pi$ have the same distribution, we obtain
\begin{align}\label{eq: second-coeffiicient-invariance}
  \E[Z_N(\beta;G)(g_J^2-1)]
  &=\E\!\left[Z_N(\beta;G_\pi)
    \bigl((\epsilon_\pi(J)g_{\pi(J)})^2-1\bigr)\right]\notag\\
  &=\E[Z_N(\beta;G)(g_{J'}^2-1)].
\end{align}}

On the other hand,  since $
  Z_N(\beta;tG)
  =Z_N(t\beta;G),
$ 
 differentiating at $t=1$ gives
 \begin{equation}\label{eq:partial-derivative-relaton}
 \sum_Ig_I\partial_{g_I}Z_N(\beta;G)
=\beta\,\partial_\beta Z_N(\beta;G).
 \end{equation}
{Therefore, substituting ${\rm He}_2(x)=x^2-1$ we obtain
\begin{align*}
	P_2Z_N(\beta)&\overset{\eqref{eq:odd-coefficient-zero}}{=}\sum_{I\in\calI_N}\frac{\E[Z_N(\beta;G)(g_I^2-1)]}{2}(g_I^2-1)\\
	&\overset{\eqref{eq: second-coeffiicient-invariance}}{=}\frac{1}{2|\calI_N|}\E\left[\sum_{I\in\calI_N}Z_N(\beta;G)(g_I^2-1)\right]\sum_{I\in\calI_N}(g_I^2-1)\\
	&=\frac{1}{2|\calI_N|}\E\left[\sum_{I\in\calI_N}g_I\partial_{g_I}Z_N(\beta;G)\right]\sum_{I\in\calI_N}(g_I^2-1)\\
	&\overset{\eqref{eq:partial-derivative-relaton}}{=}\frac{\beta\E[\partial_\beta Z_N(\beta;G)]}{2|\calI_N|}\sum_{I\in\calI_N}(g_I^2-1)\\
	&=\frac{\beta\partial_\beta\E[ Z_N(\beta;G)]}{2|\calI_N|}\sum_{I\in\calI_N}(g_I^2-1)
\end{align*}
where the third equality follows from the (one-dimensional) Gaussian integration by parts, and the last equality follows from the dominated convergence theorem because
\[\sup_{|\beta|<B}|\partial_\beta Z_N(\beta;G)|\le \|H_N(G)\|_\infty e^{B\|H_N(G)\|_\infty},\quad  \|H_N(G)\|_\infty\le \sigma_N\sum_{I\in\calI_N}|g_I|\]
for every $B>0$.} 
This proves {$P_2 Y_N(\beta)
	=A_N(\beta)Q_N$, because $P_2Y_N(\beta)=P_2Z_N(\beta)/\E Z_N(\beta)$}.
\end{proof}

\section{Direct Rademacher control of all higher chaoses} \label{sec:Rademacher-remainder}

Following the expressions of the first and second terms in the polynomial chaos expansion \eqref{eq: chaos} of $Z_N(\beta)$,
we estimate the higher-order terms in this section. Recall the definition of $Y_N(\beta)$ in \eqref{eq:Y_N} and the orthogonal projection $P_k$ defined in Section~\ref{sec:chaos-decomposition}.

\begin{proposition} \label{prop:higher-chaos-closure}
  There exists $\beta_{\mathrm{r}}=\beta_{\mathrm{r}}(q)>0$ such that, for every fixed real $\beta$ with $|\beta|<\beta_{\mathrm{r}}$,
  \begin{equation*}
    \sum_{k \geq 3}
    \E\mleft[\left|P_k Y_N(\beta)\right|^2\mright]
    =
      O\mleft(N^{-3(q/2-1)}\mright)
    =
      o\mleft(N^{2-q}\mright).
  \end{equation*}
\end{proposition}

To prove this proposition, we use the following two lemmas.

\begin{lemma} \label{lem:Z-derivative-bound}
  For any $\beta\in \mathbb R$ and for any multi-index $\alpha$, it holds 
  \begin{equation}\label{eq:Z-derivative-bound}
    \left|\partial_G^\alpha Z_N(\beta;G)\right|
    \leq (|\beta|\sigma_N)^{|\alpha|} Z_N(\beta;G).
  \end{equation}
\end{lemma}

\begin{proof}
  Let $k = |\alpha|$.
  The case that $k=0$ is immediate, so assume that $k \geq 1$ in the following. 
  We claim that, for any $(I_i)_{i=1}^{k} \subset \mathcal I_N$,
  \begin{multline}
    \partial_{g_{I_k}}\cdots\partial_{g_{I_1}}Z_N(\beta;G) \\
    =
      (-\beta\sigma_N)^k
      \sum_{\pi\in \mathfrak{S}_k}\int_{\Delta_k}
      \tr\mleft(
        e^{-t_0\beta H_N}
        \psi_{I_{\pi(1)}}
        e^{-t_1\beta H_N}
        \cdots
        \psi_{I_{\pi(k)}}
        e^{-t_k\beta H_N}
      \mright)
      \odif{t_1}\cdots\odif{t_k},
    \label{eq:iterated-Duhamel-formula}
  \end{multline}
  where
  \[
    \Delta_k =
    \left\{
      (t_j)_{j=1}^{k} \in \interval{0}{1}^k
      :
      \sum_{j=1}^{k} t_j\le 1
    \right\},
    \qquad
    t_0 = 1 - \sum_{j=1}^{k} t_j.
  \]
  By the equality that $\partial_{g_I}H_N = \sigma_N\psi_I$ and the Duhamel formula~\cite[Corollary~3.11]{HiaiPetz2014}, for any $I \in \mathcal I_N$ and any $t \geq 0$,
  \begin{align}
    \partial_{g_I}e^{-t\beta H_N}
    &=
      -t\beta\sigma_N\int_0^1
      e^{-st\beta H_N}\psi_I
      e^{-(1-s)t\beta H_N}\,\odif{s} \notag\\
    &=
      -\beta\sigma_N\int_0^t
      e^{-u\beta H_N}\psi_I
      e^{-(t-u)\beta H_N}\,\odif{u}.
    \label{eq:Duhamel-interval-splitting}
  \end{align}
  Note that  \eqref{eq:Duhamel-interval-splitting} with $t=1$ yields
  \eqref{eq:iterated-Duhamel-formula} for $k=1$.
  Suppose inductively that the formula holds at order $k - 1$.
  Then
  \begin{multline}
    \partial_{g_{I_k}}\cdots\partial_{g_{I_1}}
    Z_N(\beta;G)
    =
      (-\beta\sigma_N)^{k-1}
      \sum_{\pi\in \mathfrak{S}_{k-1}}
      \int_{\Delta_{k-1}} \odif{t_1}\cdots\odif{t_{k-1}} \\
      \times
      \partial_{g_{I_k}}
      \tr\Bigl(
        e^{-t_0\beta H_N}
        \psi_{I_{\pi(1)}}
        e^{-t_1\beta H_N}
        \cdots
        \psi_{I_{\pi(k-1)}}e^{-t_{k-1}\beta H_N}
      \Bigr)
    \label{eq:iterated-Duhamel-formula-induction}
  \end{multline}
  Fix $\pi \in \mathfrak{S}_{k-1}$ and $j\in\{1, \ldots, k\}$.
  Differentiating the factor $e^{-t_{j-1}\beta H_N}$,
  and relabeling the running variables,
  the Leibniz rule leads to
  \begin{align*}
    -\beta\sigma_N
    \int_{\Delta_k}
    \tr\mleft(
      e^{-t_0\beta H_N}
      \psi_{I_{\tilde{\pi}_{j}(1)}}
      e^{-t_1\beta H_N}
      \cdots
      \psi_{I_{\tilde{\pi}_{j}(k)}}
      e^{-t_k\beta H_N}
    \mright)
    \odif{t_1} \cdots \odif{t_k},
  \end{align*}
  where $\tilde\pi_j\in\mathfrak{S}_k$ is defined by
  \[
    \tilde{\pi}_{j}(r)
    =
    \begin{cases}
      \pi(r),   & r<j, \\
      k,        & r=j, \\
      \pi(r-1), & r>j,
    \end{cases}\qquad r\in[k].
  \]
  Since $(\pi, j) \mapsto \tilde{\pi}_{j}$ is a bijection from
  $\mathfrak{S}_{k-1}\times\{1,\ldots,k\}$ onto $\mathfrak{S}_k$, summing the terms produced by the Leibniz rule in
  \eqref{eq:iterated-Duhamel-formula-induction} proves
  \eqref{eq:iterated-Duhamel-formula}.

  To bound the integrand in \eqref{eq:iterated-Duhamel-formula}, use the normalized Schatten norms for a matrix $A$ and $p \in \interval[open right]{1}{\infty}$:
  \[
    \norm{A}_p \coloneqq \left(\tr\mleft(\sqrt{A^* A}^p\mright)\right)^{1/p}.
  \]
  Note that the case that $p=\infty$ corresponds to the operator norm $\norm{A}_\infty$.
  Let $t_0, \ldots, t_k \in \interval{0}{1}$ satisfying $\sum_{j=0}^k t_j = 1$. 
  Iterating H\"older's inequality~\cite[Theorem~6.20]{HiaiPetz2014} for any matrices $A_0, \ldots, A_k$ and $B_1, \ldots, B_k$ with the finiteness of the Schatten norm provides
  \begin{align*}
    \left|\tr(A_0 B_1 A_1 \cdots B_k A_k)\right|
    & \leq
      \norm{A_0 B_1 A_1 \cdots B_k A_k}_1\\
    & \leq
      \norm{A_0 B_1 A_1 \cdots A_{k-1} B_k}_{1/(t_0+\cdots+t_{k-1})}
      \norm{A_k}_{1/t_k} \\
    & \leq
      \norm{A_0 B_1 A_1 \cdots B_{k-1} A_{k-1}}_{1/(t_0+\cdots+t_{k-1})}
      \norm{B_k}_\infty
      \norm{A_k}_{1/t_k} \\
    & \leq
      \cdots \\
    & \leq
      \prod_{j=0}^{k} \norm{A_j}_{1/t_j}
      \prod_{i=1}^{k} \norm{B_i}_\infty.
  \end{align*}

  We apply this inequality with $A_j = e^{-t_j \beta H_N}$ and $B_i = \psi_{I_{\pi(i)}}$.
  Note that $H_N$ is Hermitian and $\norm{\psi_I}_\infty = 1$ due to $\psi_I^*\psi_I =\mathbb I$. 
  Thus, for every $\pi \in \mathfrak{S}_k$ and $(t_j)_{j=1}^{k} \in \Delta_k$,
  \begin{align}
    \MoveEqLeft
    \left|
      \tr\mleft(
        e^{-t_0\beta H_N}
        \psi_{I_{\pi(1)}}
        e^{-t_1\beta H_N}
        \cdots
        \psi_{I_{\pi(k)}}
        e^{-t_k\beta H_N}
      \mright)
    \right| \notag \\
    \leq {} &
      \prod_{j=0}^{k}
      \norm{e^{-t_j\beta H_N}}_{1/t_j}
      \prod_{i=1}^{k}
      \norm{\psi_{I_{\pi(i)}}}_\infty \notag\\
    \leq {} &
      \prod_{j=0}^{k}
      \tr\mleft(e^{-\beta H_N}\mright)^{t_j}
    =
      \tr\mleft(e^{-\beta H_N}\mright)
    =
      Z_N(\beta;G),
    \label{eq:Duhamel-integrand-Holder-bound}
  \end{align}
  where we used $\sum_{j=0}^{k} t_j=1$.
  Combining \eqref{eq:iterated-Duhamel-formula} and \eqref{eq:Duhamel-integrand-Holder-bound}, we obtain
  \begin{align*}
  \left|
    \partial_{g_{I_k}}
    \cdots
    \partial_{g_{I_1}}
    Z_N(\beta;G)
  \right|
  \leq {} &
    \left(|\beta| \sigma_N\right)^k
    \sum_{\pi\in \mathfrak{S}_k}
    \int_{\Delta_k}
    Z_N(\beta;G)
    \odif{t_1} \cdots \odif{t_k} \\
  = {} &
    \left(|\beta| \sigma_N\right)^k
    Z_N(\beta;G)
    |\mathfrak{S}_k|
    \int_{\Delta_k} 1 \odif{t_1} \cdots \odif{t_k} \\
  = {} &
    \left(|\beta| \sigma_N\right)^k
    Z_N(\beta;G),
  \end{align*}
  which proves \eqref{eq:Z-derivative-bound}.
\end{proof}

Recall that Lemma~\ref{lem: parity rule lemma} shows that only even multi-indices have nonzero Hermite coefficients.
Let $(\varepsilon_x)_{x=1}^N$ be independent Rademacher variables on an auxiliary probability space.
Write $\mathbb{P}_\varepsilon$ and $\E_\varepsilon$ for the probability measure on this space and the associated expectation, respectively.
Set
\begin{equation}\label{eq:Rademacher-chaos-definition}
  S_{N,q}
  \coloneqq
    \sum_{I\in\calI_N}\varepsilon_I,
  \qquad
  \varepsilon_I
  \coloneqq
    \prod_{x\in I}\varepsilon_x.
\end{equation}

\begin{lemma} \label{lem:Rademacher-encoding}
  For every integer $k \geq 0$,
  \begin{equation}
    \label{eq:even-sum-Rademacher-moment}
    \sum_{
      \substack{
        |\alpha|=k \\
        \alpha\ \mathrm{even}
      }
    }
    \frac{1}{\alpha!}
    =
      \frac{1}{k!}
      \E_\varepsilon\mleft[S_{N,q}^k\mright].
  \end{equation}
  Moreover, there are constants $c_q > 0$ and $C_q < \infty$ depending only on $q$ such that, for every $N$ and every $0 \leq t \leq c_q N^{-(q-1)}$,
  \begin{equation}
    \label{eq:Rademacher-exponential-moment}
    \E_\varepsilon\mleft[
      |S_{N,q}|^3
      e^{t |S_{N,q}|}
    \mright]
    \leq
      C_q N^{3q/2}.
  \end{equation}
\end{lemma}

\begin{proof}
  We first show \eqref{eq:even-sum-Rademacher-moment}.
  Starting from the definition \eqref{eq:Rademacher-chaos-definition}
  and applying the multinomial formula provides
  \begin{equation}\label{eq:Rademacher-multinomial-expansion}
    \frac{1}{k!}\E_\varepsilon\mleft[S_{N,q}^k\mright]
    =\frac{1}{k!}\E_\varepsilon\mleft[
      \sum_{|\alpha|=k}k!
      \prod_{I\in\calI_N}\frac{\varepsilon_I^{\alpha_I}}{\alpha_I!}
    \mright]
    =
    \sum_{|\alpha|=k}\frac{1}{\alpha!}
    \E_\varepsilon\mleft[
      \prod_{I\in\calI_N}\varepsilon_I^{\alpha_I}
    \mright].
  \end{equation}
  For a fixed $\alpha$, the product inside the expectation is
  \begin{equation*}
    \prod_{x=1}^N
    \varepsilon_x^{\sum_{I: x\in I}\alpha_I}.
  \end{equation*}
  By independence of the Rademacher variables and the definition of an even
  multi-index in Section~\ref{sec:chaos-decomposition}, its expectation is one if
  $\alpha$ is even and zero otherwise.  Substituting this into
  \eqref{eq:Rademacher-multinomial-expansion} proves
  \eqref{eq:even-sum-Rademacher-moment}.

  We next bound $S_{N,q}$.  For every integer $N\geq1$ and $m\geq0$, let $e_m^{(N)}$ denote the $m$th elementary symmetric polynomial in $N$ variables:
  \begin{equation*}
    e_m^{(N)}(z_1,\ldots,z_N)
    \coloneqq
      \sum_{
        \substack{
          I\subset[N]\\
          |I|=m
        }
      }
      \prod_{x\in I}
      z_x,
    \qquad
    e_0^{(N)}(z_1,\ldots,z_N)
    \coloneqq
      1.
  \end{equation*}
  For the rest of the proof, we abbreviate $e_m \coloneqq e_m^{(N)}(\varepsilon_1,\ldots,\varepsilon_N)$.
  Then, $e_1=\sum_{x=1}^N\varepsilon_x$ and \eqref{eq:Rademacher-chaos-definition} imply that $S_{N,q}=e_q$.
  Newton's identities (see, e.g., \cite[Chapter I, \S2]{Macdonald1995}), state that
  \begin{equation}\label{eq:Newton-identities}
    m e_m
    =\sum_{j=1}^m(-1)^{j-1}e_{m-j}p_j,
  \end{equation}
  where $p_j \coloneqq \sum_{x=1}^N\varepsilon_x^j$.  Since
  $\varepsilon_x\in\{\pm1\}$, we have
  \begin{equation*}
    p_j
    =
      \begin{cases}
        e_1, & j\text{ odd},  \\
        N,   & j\text{ even}.
      \end{cases}
  \end{equation*}
  We now prove by induction on $m$ that there is a constant $C_m < \infty$ such that
  \begin{equation*}
    |e_m|
    \le C_m\left(N^{m/2}+|e_1|^m\right).
  \end{equation*}
  Suppose that the estimate holds for  orders less than $m$.
  Set $X_N \coloneqq e_1 / \sqrt N$.
  For every integer $j\ge1$,
  \[
    |p_j| \le N^{j/2} \left(1 + |X_N|^j\right),
  \]
  hence, by the induction hypothesis and putting $C_0 = 1$,
  \begin{align*}
    |e_{m-j}p_j|
    &\le
      C_{m-j}N^{m/2}
      (1+|X_N|^{m-j})(1+|X_N|^j)\\
    &\le
      2C_{m-j}N^{m/2}(1+|X_N|^m),
  \end{align*}
  where we used $1+|X|^m-|X|^{m-j}-|X|^j = (1-|X|^{m-j})(1-|X|^j) \ge 0$.
  Substituting this into \eqref{eq:Newton-identities} proves the induction step.
  Thus, taking $m=q$ and using \eqref{eq:Rademacher-chaos-definition}, we obtain
  \begin{equation}\label{eq:Rademacher-chaos-deterministic-bound}
    |S_{N,q}|
    \leq C_q\left(N^{q/2}+|e_1|^q\right).
  \end{equation}

  Since $q\geq4$ and $|X_N|\leq\sqrt N$,
  \eqref{eq:Rademacher-chaos-deterministic-bound} implies that, for $0\leq t\leq c_qN^{-(q-1)}$,
  \begin{equation}
    t|S_{N,q}|
    \leq
      C_q c_q N^{-(q/2-1)} \left(1 + |X_N|^q\right)
    \leq
      C_q c_q \left(1 + |X_N|^2\right).
    \label{eq:Rademacher-exponent-reduction}
  \end{equation}
  Indeed,
  \begin{equation*}
    N^{-(q/2-1)}|X_N|^q
    =|X_N|^2
    \left(\frac{|X_N|}{\sqrt N}\right)^{q-2}
    \leq |X_N|^2.
  \end{equation*}
  Choose $c_q$ small enough so that $C_q c_q\leq 1/4$.  Hoeffding's inequality~\cite[Theorem~2]{Hoeffding1963} provides
  \begin{equation}\label{eq:Hoeffding-Rademacher-sum}
    \mathbb P_\varepsilon(|X_N|\geq u)
    \leq
      2e^{-u^2/2}
    \qquad
    (u\geq0).
  \end{equation}
  Consequently, integrating the tail bound in \eqref{eq:Hoeffding-Rademacher-sum} yields
  \begin{equation}\label{eq:uniform-Rademacher-integrability}
    \sup_{N\geq1}
    \E_\varepsilon\mleft[
      \left(1+|X_N|^{3q}\right)
      e^{|X_N|^2/4}
    \mright]
    <
    \infty.
  \end{equation}
  Finally, \eqref{eq:Rademacher-chaos-deterministic-bound}, \eqref{eq:Rademacher-exponent-reduction} and \eqref{eq:uniform-Rademacher-integrability} imply that
  \begin{align*}
    \E_\varepsilon\mleft[
      |S_{N,q}|^3e^{t|S_{N,q}|}
    \mright]
    &\leq
    C_qN^{3q/2}
    \E_\varepsilon\mleft[
      \left(1+|X_N|^{3q}\right)
      e^{(1+|X_N|^2)/4}
    \mright] \\
    &\leq C_qN^{3q/2},
  \end{align*}
  which proves \eqref{eq:Rademacher-exponential-moment}.
\end{proof}

\begin{proof}[Proof of Proposition~\ref{prop:higher-chaos-closure}]
  For any $k = |\alpha| \geq 3$, the polynomial chaos expansion \eqref{eq: chaos} and its orthogonality relation lead to
  \[
    \E\mleft[\left|P_k Y_N(\beta)\right|^2\mright]
    =
      \frac{1}{\E[Z_N(\beta)]^2}
      \sum_{|\alpha|=k}\frac{1}{\alpha!}
      \left|
        \E\mleft[Z_N(\beta) H_\alpha(G)\mright]
      \right|^2.
  \]
  By Lemmas~\ref{lem: parity rule lemma} and~\ref{lem:Gaussian-Hermite-IBP},
  \[
    \E\mleft[Z_N(\beta)H_\alpha(G)\mright]
    =
      \begin{cases}
        0,                                & \alpha \text{ is odd},  \\
        \E[\partial_G^\alpha Z_N(\beta)], & \alpha \text{ is even}.
      \end{cases}
  \]
  Moreover, Lemma~\ref{lem:Z-derivative-bound} provides,
  \[
    \bigl|\E[\partial_G^\alpha Z_N(\beta)]\bigr|
    \leq
      \E\bigl[\left|\partial_G^\alpha Z_N(\beta)\right|\bigr]
    \leq
      (|\beta|\sigma_N)^{k}
      \E\mleft[Z_N(\beta)\mright].
  \]
  Since $\E[Z_N(\beta)] > 0$, it follows that
  \begin{equation}
    \E\mleft[\left|P_k Y_N(\beta)\right|^2\mright]
    =
      \frac{1}{\E[Z_N(\beta)]^2}
      \sum_{
        \substack{
          |\alpha| = k \\
          \alpha\ \mathrm{even}
        }
      }
      \frac{1}{\alpha!}
      \bigl|\E\mleft[\partial_G^\alpha Z_N(\beta)\mright]\bigr|^2
    \leq (|\beta| \sigma_N)^{2 k}
      \sum_{
        \substack{
          |\alpha| = k \\
          \alpha\ \mathrm{even}
        }
      }
      \frac{1}{\alpha!}.
    \label{eq:Hermite-Rademacher-bound}
  \end{equation}
  By Lemma~\ref{lem:Rademacher-encoding}, summing \eqref{eq:Hermite-Rademacher-bound} over all $k \geq 3$ yields
  \begin{align}
    \sum_{k\geq3}\E\mleft[\left|P_k Y_N(\beta)\right|^2\mright]
    & \leq
      \sum_{k\geq3}
      \frac{(|\beta| \sigma_N)^{2 k}}{k!}
      \E_\varepsilon\mleft[S_{N,q}^k\mright] \notag\\
    & \leq
      \E_\varepsilon\mleft[
        \sum_{k\geq3}
        \frac{(\beta^2 \sigma_N^2 |S_{N,q}|)^k}{k!}
      \mright] \notag\\
    & \leq
      \frac{(|\beta| \sigma_N)^{6}}{6}
      \E_\varepsilon\mleft[
        |S_{N,q}|^3
        e^{\beta^2 \sigma_N^2 |S_{N,q}|}
      \mright],
    \label{eq:all-higher-chaoses-at-once}
  \end{align}
  where the second inequality uses $|\E_\varepsilon[S_{N,q}^k]| \leq \E_\varepsilon[|S_{N,q}|^k]$,
  and the last inequality uses $\sum_{k\geq3}u^k / k! \leq u^3e^u/6$ for $u\geq0$.
  Choose
  \begin{equation*}
    \beta_{\mathrm{r}}
    \coloneqq
      \sqrt{\frac{c_q}{(q-1)!}},
  \end{equation*}
  where $c_q$ is the constant in Lemma~\ref{lem:Rademacher-encoding}.
  If $|\beta| < \beta_{\mathrm{r}}$, then $\beta^2\sigma_N^2 = \beta^2(q-1)!N^{-(q-1)} \leq c_q N^{-(q-1)}$.
  Therefore the second claim of Lemma~\ref{lem:Rademacher-encoding}, together with \eqref{eq:all-higher-chaoses-at-once}, implies
  \begin{align*}
    \sum_{k\geq3}
    \E\mleft[
      \left|P_k Y_N(\beta)\right|^2
    \mright]
    &\leq
      C_{q, \beta}
      N^{-3(q-1)}
      N^{3q/2}
    =
      C_{q, \beta}
      N^{-3(q/2-1)}.
  \end{align*}
  Since $q \geq 4$, we have $3 (q / 2 - 1) > q - 2$, which completes the proof.
\end{proof}

\section{Proof of the main results}
\label{sec:main-results}

To identify the deterministic coefficient $A_N(\beta)$ in
\eqref{eq:finite-N-FTW-coefficient}, we use the annealed free energy limit
from \cite{He2026}, which holds at every fixed non-negative $\beta$. 
The symmetry of Gaussian distributions extends the claim to negative \(\beta\).

The following is a standard convergence theorem for convex functions, {so we omit its proof}; see
\cite[Theorems~10.8 and~25.7]{Rockafellar1970}.

\begin{lemma} \label{lem:convex-derivatives}
Let $h\in C^1(\mathbb R )$. Suppose that $h_N:\mathbb R \to\mathbb R$ are
 differentiable, convex functions, and that $h_N\to h$ pointwise. Then locally uniformly on $\mathbb R $, 
\[
    h_N\to h
    \qquad\text{and}\quad
    h_N'\to h'.
\]

\end{lemma}
{The next proposition establishes convergence of the annealed free energy and its derivative.}
\begin{proposition} \label{prop:annealed-free-energy}
    There exists an even convex function $\ell_q\in C^1(\mathbb R)$ such that,
    locally uniformly for $\beta\in\mathbb R$,
    \begin{align}
        \frac1N\log \E[Z_N(\beta)]
        &\longrightarrow \ell_q(\beta),
        \label{eq:annealed-free-energy-limit}\\
        \frac1N\partial_\beta\log \E[Z_N(\beta)]
        &\longrightarrow \ell_q'(\beta).
        \label{eq:annealed-derivative-limit}
    \end{align}
    Moreover, $\beta\ell_q'(\beta)>0$ for every $\beta\neq0$.
\end{proposition}

\begin{proof}
    Set
    \begin{equation*}
        \widetilde{H}_N:=2^{-q/2}H_N,
        \qquad
       {\widetilde{p}_{N,q}(s)}:=\frac1N\log\E[\tr e^{s\widetilde{H}_N}].
    \end{equation*}
    Since $G$ and $-G$ have the
    same law and $H_N(-G)=-H_N(G)$, for every $\beta\in\mathbb R$,
    \begin{equation} \label{eq:He-normalization-conversion}
        \frac1N\log\E[Z_N(\beta)]
        =\frac1N\log\E[\tr e^{|\beta|H_N}]
        =\widetilde{p}_{N,q}(2^{q/2}|\beta|).
    \end{equation}
    Writing $p_q:=p_{{\rm SD},q}$ for the Schwinger--Dyson pressure, \cite[Theorem~1.1 and Proposition~3.3]{He2026} give
    \begin{equation*}
        \widetilde{p}_{N,q}(s)\longrightarrow p_q(s)
        \quad (s>0),
        \qquad
        p_q\in C^1((0,\infty)),
        \qquad
        \lim_{s\downarrow0}p_q(s)=0.
    \end{equation*}
    Define
    \begin{equation*}
        \ell_q(\beta):=
        \begin{cases}
            p_q(2^{q/2}|\beta|),&\beta\neq0,\\
            0,&\beta=0.
        \end{cases}
    \end{equation*}
    Then $\ell_q$ is even and continuous on $\mathbb R$. By \eqref{eq:He-normalization-conversion} and  $Z_N(0)=1$, we have pointwise convergence in \eqref{eq:annealed-free-energy-limit}.

    By the derivative formula in \cite[Proposition 3.3]{He2026}, for any $s>0$, 
    \begin{equation} \label{eq:He-pressure-derivative}
        p_q'(s)=\frac{s}{q}\int_0^1G_{s,q}(u)^q\,\mathrm du.
    \end{equation}
    By the positive Lehmann representation and
    \cite[Theorem~2.5]{He2026}, there is a probability measure $\nu_s$
    on $[0,\infty)$ such that
    \begin{equation*}
        G_{s,q}(u)
        =\int_0^\infty
        \frac{\cosh(x(1/2-u))}{2\cosh(x/2)}\,\nu_s(\mathrm dx).
    \end{equation*}
    For $\beta>0$, put
    \begin{equation*}
        K_\beta(u):=2G_{2^{q/2}\beta,q}(u).
    \end{equation*}
    The chain rule and
    \eqref{eq:He-pressure-derivative} yield
    \begin{equation} \label{eq:ell-kernel-derivative}
        \ell_q'(\beta)
        =2^{q/2}p_q'(2^{q/2}\beta)
        =\frac{2^q\beta}{q}\int_0^1G_{2^{q/2}\beta,q}(u)^q\,\mathrm du
        =\frac\beta q\int_0^1K_\beta(u)^q\,\mathrm du.
    \end{equation}
    Moreover,
    \begin{equation} \label{eq:ell-kernel-representation}
        K_\beta(u)
        =\int_0^\infty
        \frac{\cosh(x(1/2-u))}{\cosh(x/2)}\,\nu_{2^{q/2}\beta}(\mathrm dx).
    \end{equation}

    For $x\geq0$ and $u\in[0,1]$,
    \begin{equation*}
        0<\frac{\cosh(x(1/2-u))}{\cosh(x/2)}\leq1.
    \end{equation*}
    Since $\nu_{2^{q/2}\beta}$ is a probability measure,
    \eqref{eq:ell-kernel-representation} gives $0<K_\beta(u)\leq1$.
    Hence \eqref{eq:ell-kernel-derivative} implies
    \begin{equation*}
        0<\ell_q'(\beta)\leq\frac\beta q, 
        \qquad \beta>0.
    \end{equation*}
    Integrating over $[0,\beta]$ with $\ell_q(0)=0$ 
    {gives} $0\leq\ell_q(\beta)\leq\beta^2/(2q)$ for $\beta>0$.
    Since $\ell_q$ is even, this bound holds for every $\beta\in\mathbb R$, thus 
    \begin{equation*}
        \ell_q'(0)= \lim_{\beta\to0}
        \frac{\ell_q(\beta)-\ell_q(0)}{\beta}=0.
    \end{equation*}

On $\mathbb R\setminus\{0\}$, the derivative
    $\ell_q'$ is continuous by $p_q\in C^1((0,\infty))$ and is odd
    by evenness of $\ell_q$. Therefore
    $|\ell_q'(\beta)|\leq|\beta|/q$, which gives continuity at zero
    and hence $\ell_q\in C^1(\mathbb R)$.
    Strict positivity for $\beta>0$ and oddness also give
    $\beta\ell_q'(\beta)>0$ for every $\beta\neq0$.

It remains to prove convexity. 
    Set
    \begin{equation*}
        F_N(\beta):=\E[Z_N(\beta)],
        \qquad
        h_N(\beta):=\frac1N\log F_N(\beta),
        \qquad
        \overline\nu_N(B):=\E[\tr\mathbf1_B(H_N)].
    \end{equation*}
    Then $\overline\nu_N$ is a probability measure and
    \begin{equation*}
        F_N(\beta)=\int_{\mathbb R}e^{-\beta x}\,
        \overline\nu_N(\mathrm dx).
    \end{equation*}
 In particular,
    $h_N$ is finite and twice continuously differentiable.  With the
    probability measure
    \begin{equation*}
        \overline\nu_{N,\beta}(\mathrm dx)
        :=\frac{e^{-\beta x}}{F_N(\beta)}\,
        \overline\nu_N(\mathrm dx),
    \end{equation*}
    we have
    \begin{equation*}
        h_N''(\beta)
        =\frac1N\left(
            \int_{\mathbb R}x^2\,\overline\nu_{N,\beta}(\mathrm dx)
            -\left(\int_{\mathbb R}x\,
            \overline\nu_{N,\beta}(\mathrm dx)\right)^2
        \right)\geq0.
    \end{equation*}
    Thus, $h_N$ is convex. 
    Applying Lemma~\ref{lem:convex-derivatives} proves
    both \eqref{eq:annealed-free-energy-limit} and
    \eqref{eq:annealed-derivative-limit} locally uniformly on $\mathbb R$.
\end{proof}

By \eqref{eq:finite-N-FTW-coefficient} and
\eqref{eq:annealed-derivative-limit}, locally uniformly for
$\beta\in\mathbb R$,
\begin{equation*}
  \frac{A_N(\beta)}N\longrightarrow
  \frac\beta2\ell_q'(\beta).
\end{equation*}

We are now in a position to prove Theorems~\ref{thm:partition-clt} and~\ref{thm:free-energy-clt}.

\begin{proof}[Proof of Theorem~\ref{thm:partition-clt}]
Note that $\beta\ell_q'(\beta)>0$ from Proposition~\ref{prop:annealed-free-energy}. We take $\beta_q=\beta_{\mathrm r}$ from Proposition~\ref{prop:higher-chaos-closure} and fix $0<|\beta|<\beta_q$. Let
\[
  U_N:=\frac1{\sqrt{2|\mathcal I_N|}}\sum_{I\in\calI_N}(g_I^2-1).
\]
The classical i.i.d.\ central limit theorem gives
$U_N\Rightarrow\N(0,1)$.
By {Proposition~\ref{prop:exact-second-chaos}},
\begin{equation}\label{eq:Y-decomposition-final}
  Y_N(\beta)-1
  =c_N(\beta)U_N+\mathcal R_N(\beta),
  \qquad
  c_N(\beta)
  :=\frac{\beta\,\partial_\beta\log \E[Z_N(\beta)]}{\sqrt{2|\mathcal I_N|}},
\end{equation}
where $\mathcal R_N:=\sum_{r\geq3}P_r Y_N$ is orthogonal to $U_N$.
Since $|\mathcal I_N|=\binom Nq\sim N^q/q!$, \eqref{eq:annealed-derivative-limit} gives
\[
  N^{q/2-1}c_N(\beta)
  \longrightarrow
  \sqrt{\frac{q!}{2}}\,\beta\ell_q'(\beta).
\]
Proposition~\ref{prop:higher-chaos-closure} gives
$N^{q/2-1}\mathcal R_N\to0$ in $L^2$.  This proves
\eqref{eq:scaled-relative-clt}.  Orthogonality in
\eqref{eq:Y-decomposition-final} also gives
\[
  \Var(Y_N)=c_N(\beta)^2+\norm{\mathcal R_N(\beta)}_2^2.
\]
Thus $N^{q-2}\Var(Y_N)\longrightarrow(q!/2)\beta^2(\ell_q'(\beta))^2>0$.
Together with \eqref{eq:scaled-relative-clt} and the identity
$\Var(Y_N)=\Var(Z_N)/(\E[Z_N])^2$, Slutsky's theorem gives
\eqref{eq:standardized-clt}.
\end{proof}

\begin{proof}[Proof of {Theorem~\ref{thm:free-energy-clt}}]
Assume $\beta\neq0$.  Set
\[
a:=\frac q2-1,
\qquad
W_N:=\log Y_N(\beta).
\]
Then
\begin{align*}
N^a\bigl(\log Z_N-\E[\log Z_N]\bigr)
&=
N^a\bigl(W_N-(Y_N-1)\bigr)
+N^a(Y_N-1)
-N^a\E[W_N] .
\end{align*}

By Theorem~\ref{thm:partition-clt},
\[
N^a(Y_N-1)
\Longrightarrow
\sqrt{\frac{q!}{2}}\,
\beta\ell_q'(\beta)\,\xi,
\qquad
\xi\sim\mathcal N(0,1).
\]
Moreover, $Y_N-1\to0$ in probability and
$N^a(Y_N-1)$ is tight.  Since whenever $|x|\le1/2$, 
\[
|\log(1+x)-x|\le x^2,
\]
it follows that
\[
N^a\bigl(W_N-(Y_N-1)\bigr)
\xrightarrow{\mathrm p}0.
\]
Indeed, on $\{|Y_N-1|\leq1/2\}$ the absolute value is bounded by
\(
|N^a(Y_N-1)|\,|Y_N-1|
\), which converges to zero in probability by tightness of the first factor
and convergence of the second factor to zero.  The complementary event has
probability tending to zero.
Thus it remains to prove
\[
\E[W_N]=o(N^{-a}).
\]

Write, as in \eqref{eq:Y-decomposition-final},
\[
Y_N-1=c_NU_N+\mathcal R_N,
\]
where
\[
c_N\coloneqq \frac{\beta\,\partial_\beta\log \E[Z_N(\beta)]}{\sqrt{2|\mathcal I_N|}}= O(N^{-a}),\quad U_N\coloneqq \frac1{\sqrt{2|\mathcal I_N|}}\sum_{I\in\mathcal I_N}(g_I^2-1)
,\quad
\|\mathcal R_N\|_2^2=O(N^{-3a})
\]
by Proposition~\ref{prop:higher-chaos-closure}. Orthogonality therefore
gives
\[
\E[(Y_N-1)^2]=c_N^2 \E[U_N^2] + \E[\mathcal R_N^2]= O(N^{-2a}).
\]
Since $\E[(g_I^2-1)^2]=2$ and
$\E[(g_I^2-1)^4]=60$, independence gives
\[
\E[U_N^4]=3+\frac{12}{|\mathcal I_N|}\leq15.
\]  
Hence, we estimate
\begin{align*}
\mathbb P(Y_N<1/2)
&\le
\mathbb P(|c_NU_N|>1/4)
+
\mathbb P(|\mathcal R_N|>1/4)\\
&\le
C c_N^4\E[U_N^4]
+
C\|\mathcal R_N\|_2^2\\
&=O(N^{-3a}).
\end{align*}

We claim that 
\[
Z_N(\beta;G)\ge1.
\]
To this end, fix \(I\in\mathcal I_N\) and choose \(x\in I\). 
The unitary $U_x$ {constructed in the proof of}  Lemma~\ref{lem: parity rule lemma} satisfies
\[
U_x\psi_IU_x^*=-\psi_I,
\]
and hence $\tr\psi_I=0$.  Therefore
\(
\tr H_N(G)=0.
\) 
Jensen's inequality applied to the eigenvalues of $H_N(G)$ now gives
\[
Z_N(\beta;G)
=
\tr e^{-\beta H_N(G)}
\ge
e^{-\beta\tr H_N(G)}
=1.
\]

Define
\[
h(y):=y-1-\log y,
\qquad y>0.
\]
Since $\E[Y_N]=1$,
\[
-\E[W_N]
=
\E[h(Y_N)].
\]
On $\{Y_N\ge1/2\}$,
\[
h(Y_N)\le (Y_N-1)^2,
\]
and therefore
\[
\E\!\left[
h(Y_N)\mathbf 1_{\{Y_N\ge1/2\}}
\right]
=O(N^{-2a}).
\]

On the other hand, if $Y_N<1/2$, then
\[
h(Y_N)\le-\log Y_N.
\]
Since $Z_N\ge1$,
\[
-\log Y_N
=
\log\E[Z_N]-\log Z_N
\le
\log\E[Z_N].
\]
By Proposition~\ref{prop:annealed-free-energy},
\(
\log\E[Z_N]=O(N),
\) 
we have
\begin{align*}
\E\!\left[
h(Y_N)\mathbf 1_{\{Y_N<1/2\}}
\right]
&\le
O(N)\,\mathbb P(Y_N<1/2)=
O(N^{1-3a})
=O(N^{-2a}),
\end{align*}
where the last estimate uses $a=q/2-1\geq1$.
Consequently,
\[
0\le-\E[W_N]
=
O(N^{-2a})
=o(N^{-a}).
\]

Combining the three terms above and applying {Slutsky's lemma
\cite[Theorem~3.1]{Billingsley1999}} yields
\[
N^{q/2-1}
\bigl(
\log Z_N(\beta)-\E[\log Z_N(\beta)]
\bigr)
\Longrightarrow
\mathcal N\left(
0,\frac{q!}{2}\beta^2(\ell_q'(\beta))^2
\right).
\]
\end{proof}

\appendix
\section*{Lean formalization.}
We have formalized the proofs of our main results in Lean 4, conditional
on Proposition~\ref{prop:annealed-free-energy}. Proposition~\ref{prop:annealed-free-energy} is
treated as an explicit assumption and has not
been formalized. 
The Lean code is available at
\url{https://github.com/TokyoProbabilityGroup/SYK_CLT}.

\section*{Acknowledgements}
  The first author was supported by JSPS KAKENHI Grant Number JP26K17003.
 The second author was supported by JSPS KAKENHI Grant Numbers JP24K16937 and JP25K00911. 
 The third author was supported by JSPS KAKENHI Grant Number JP26KJ0336. 
  The fourth author was supported by JSPS KAKENHI Grant Numbers JP24KJ0022 and JP26K17012.

\end{document}